\documentclass[11pt, a4paper]{article}
\usepackage{amssymb}
\usepackage{amsmath}
\usepackage{amsthm}
\usepackage{rotate}
\usepackage{graphicx}

\theoremstyle{definition}
\newtheorem{definition}{Definition}[section]
\theoremstyle{plain}
\newtheorem{theorem}[definition]{Theorem}
\newtheorem{lemma}[definition]{Lemma}
\newtheorem{proposition}[definition]{Proposition}
\newtheorem{corollary}[definition]{Corollary}
\newtheorem{example}[definition]{Example}

\begin{document}

\title{\bf Intuitionistic fuzzy $H_v$-submodules}
\author{\bf Bijan Davvaz$^{\rm a}$, Wies{\l}aw A. Dudek$^{\rm b},$
             Young Bae Jun$^{\rm c,*}$ \\[2mm]
 {\small \it $^{\rm a}$ Department of Mathematics, Yazd University, Yazd, Iran}\\
 {\small \it $^{\rm b}$ Institute of Mathematics, Technical University,} \\
 {\small \it Wybrze\.ze Wyspia\'nskiego 27, 50-370 Wroc{\l}aw, Poland}\\
 {\small \it $^{\rm c}$ Department of Mathematics Educations,
               Gyeongsang National University,} \\
 {\small \it   Chinju 660-701, Korea}}

\date{}
\maketitle

\begin{abstract}
 After the introduction of fuzzy sets by Zadeh, there have been a
 number of generalizations of this fundamental concept. The notion of
intuitionistic fuzzy sets introduced by Atanassov is one among
them. In this paper, we apply the concept of an intuitionistic
fuzzy set to $H_v$-modules. The notion of an intuitionistic fuzzy
$H_v$-submodule of an $H_v$-module is introduced, and some related
properties are investigated. Characterizations of intuitionistic
fuzzy $H_v$-submodules are given. 
\end{abstract}

  \footnote{{\it E-mail address:}
        davvaz\symbol{64}yazduni.ac.ir (B. Davvaz),
        dudek\symbol{64}im.pwr.wroc.pl (W. A. Dudek),
        ybjun\symbol{64}gsnu.ac.kr (Y. B. Jun)}

\section{Introduction}

The concept of hyperstructure was introduced  in 1934 by Marty
\cite{14} at the 8th congress of Scandinavian Mathematicians.
Hyperstructures have many applications to several branches of both
pure and applied sciences \cite{4, 5}. Vougiouklis \cite{19}
introduced a new class of hyperstructures so-called
$H_v$-structure, and Davvaz \cite{9} surveyed the theory of
$H_v$-structures.  After the introduction of fuzzy sets by Zadeh
\cite{21}, there have been a number of generalizations of this
fundamental concept. The notion of intuitionistic fuzzy sets
introduced by Atanassov \cite{1} is one among them.
 In \cite{3}, Biswas applied the concept of intuitionistic fuzzy sets
to the theory of groups and studied intuitionistic fuzzy subgroups
of a group. In \cite{12}, Kim, Dudek and Jun introduced the notion
of intuitionistic fuzzy subquasigroups of a quasigroup. Also in
\cite{13}, Kim and Jun introduced the concept of intuitionistic
fuzzy ideals of a semigroup. Recently, Dudek, Davvaz and Jun
\cite{11} considered the intuitionistic fuzzification of the
concept of sub-hyperquasigroups in a hyperquasigroup and
investigated some properties of such hyperquasigroups.
 In this paper, we apply the concept of
intuitionistic fuzzy sets to $H_v$-modules. We introduce the
notion of intuitionistic fuzzy $H_v$-submodules of an $H_v$-module
and investigate some related properties. We give characterizations
of intuitionistic fuzzy $H_v$-submodules.

\section{Fuzzy sets and intuitionistic fuzzy sets} 

The concept of a fuzzy set in a non-empty set was introduced by
Zadeh \cite{21} in 1965.

Let $X$ be a non-empty set. A mapping $\mu :X\longrightarrow
[0,1]$ is called a {\it fuzzy set} in $X$. The {\it complement} of
$\mu$, denoted by $\mu^c$, is the fuzzy set in $X$ given by
$\mu^c(x)=1-\mu (x)$ for all $x\in X$.

\begin{definition}
Let $f$ be a mapping from a set $X$ to a set $Y$. Let $\mu$ be a
fuzzy set in $X$ and $\lambda$ be a fuzzy set in $Y$. Then the
{\it inverse image} $f^{-1}(\lambda )$ of $\lambda$ is a fuzzy set
in $X$ defined by
\[ f^{-1}(\lambda )(x)=\lambda (f(x)) \
{\rm for \ all} \ x \in X.
\]
The {\it image} $f(\mu )$ of $\mu$ is the fuzzy set in $Y$ defined
by
\[ f(\mu )(y)=\left \{ \begin{array}{cl} \displaystyle
   \sup_{x \in f^{-1}(y)}\mu (x) & {\rm if } \ f^{-1}(y)\ne \emptyset, \\
   0 & {\rm otherwise,}
\end{array}\right. \]
for all $y \in Y$. We have always
\[
f(f^{-1}(\lambda ))\leq \lambda \ \ {\rm and} \ \ \mu \leq
f^{-1}(f(\mu )).
\]
\end{definition}

Rosenfeld \cite{16} applied the concept of fuzzy sets to the
theory of groups and defined the concept of fuzzy subgroups of a
group.  The concept of fuzzy modules was introduced by Negoita and
Ralescu in \cite{15}.

\begin{definition} (cf. Negoita and Ralescu \cite{15}). Let $M$ be a
module over a ring $R$. A fuzzy set $\mu$ in $M$ is called a {\em
fuzzy submodule} of $M$ if for every $x,y \in M$ and $r \in R$ the
following conditions are satisfied:
\begin{itemize}
 \item[(1)]$\mu (0)=1$,
 \item[(2)] $\min \{ \mu (x), \mu (y)\}\leq \mu (x-y)$ for all $x,y \in M$,
 \item[(3)] $\mu (x)\leq \mu (rx)$ for all $x \in M$ and $r\in R$.
\end{itemize}
\end{definition}

\begin{definition} An {\em intuitionistic
fuzzy set} $A$ in a non-empty set $X$ is an object having the form
\[
A=\{(x,\mu_A(x),\lambda_A(x)) \ | \ x\in X\},
\]
where the functions $\mu_A:X\longrightarrow [0,1]$ and
$\lambda_A:X\longrightarrow [0,1]$ denote the degree of membership
(namely $\mu_A(x)$) and the degree of nonmembership (namely
$\lambda_A(x)$) of each element $x\in X$ to the set $A$
respectively, and $0\leq \mu_A(x)+\lambda_A(x)\leq  1$ for all $x
\in X$. For the sake of simplicity, we shall use the symbol
$A=(\mu_A,\lambda_A)$ for the intuitionistic fuzzy set
$A=\{(x,\mu_A(x),\lambda_A(x)) \ | \ x\in X\}$.
\end{definition}

\begin{definition} Let $A=(\mu_A,\lambda_A)$ and $B=(\mu_B,\lambda_B)$ be
intuitionistic fuzzy sets in $X.$ Then
\begin{itemize}
 \item[(1)] $ A\subseteq B$ \ iff \ $\mu_A(x)\leq\mu_B(x)$ \ and \
$\lambda_A(x)\geq\lambda_B(x)$ \ for all $x\in X$,
 \item[(2)] $A^c=\{(x,\lambda_A(x),\mu_A (x)) \ | \ x\in X\}$,
 \item[(3)] $A\cap B=\{(x,\min\{\mu_A(x),\mu_B(x)\},\max\{\lambda_A(x),\lambda_B(x)\})
\ | \ x\in X\}$,
 \item[(4)] $A\cup B=\{(x,\max\{\mu_A(x),\mu_B(x)\},\min\{\lambda_A(x),\lambda_B(x)\})
\ | \ x\in X\}$,
 \item[(5)] $\Box A=\{(x,\mu_A(x),\mu^c_A(x)) \ | \ x\in X\}$,
 \item[(6)] $\Diamond A= \{(x,\lambda^c_A(x),\lambda_A(x)) \ | \ x\in X\}$.
\end{itemize}
\end{definition}

Now, we define an intuitionistic fuzzy submodule of a module.

\begin{definition} Let $M$ be a module over a ring $R$. An
intuitionistic fuzzy set $A=(\mu_A,\lambda_A)$ in $M$ is called an
{\em intuitionistic fuzzy submodule} of $M$ if
\begin{itemize}
 \item[(1)] \ $\mu_A (0)=1$,
 \item[(2)] \ $\min\{\mu_A(x),\mu_A(y)\}\leq\mu_A(x-y)$ \ \ for all $x,y\in M$,
 \item[(3)] \ $\mu_A(x)\leq\mu_A(r\cdot x)$\ \ for all $x\in M$ and
$r \in R$,
 \item[(4)] \ $\lambda_A(0)=0$,
 \item[(5)] \ $\lambda_A(x-y)\leq\max\{\lambda_A(x),\lambda_A(y)\}$ \ \ for all
$x,y \in M$,
 \item[(6)] \ $\lambda_A(r\cdot x)\leq\lambda_A(x)$ \ \
\ \ for all $x\in M$ and $r\in R$.
\end{itemize}
\end{definition}

\section{$H_v$-structures} 

A {\em hyperstructure} is a non-empty set $H$ together with a map
$*:H\times H\to {\cal P}^*(H)$ which is called {\em
hyperoperation}, where ${\cal P}^*(H)$ denotes the set of all
non-empty subsets of $H$. The image of pair $(x,y)$ is denoted by
$x*y$. If $x\in H$ and $A,B\subseteq H$, then by $A*B$, $A*x$ and
$x*B$ we mean
\[
A*B=\displaystyle\bigcup_{a\in A, b\in B}\!\!\! a*b, \ \ A*x=A*\{
x\} \ \ {\rm and } \ \ x*B=\{x\}*B,
\]
respectively. A hyperstructure $(H,*)$ is called an {\em
$H_v$-semigroup} if
\[
(x*(y*z))\cap ((x*y)*z)\neq\emptyset \ \ {\rm for \ all } \
x,y,z\in H.
\]

\begin{definition} \label{def11} An {\em $H_v$-ring} is
a system $(R,+,\cdot )$ with two hyperoperations satisfying the
following axioms:
\begin{itemize}
\item[(i)] \ $(R,+)$ is an {\em $H_v$-group}, i.e.,
\[
\begin{array}{l}
((x+y)+z )\cap ( x+(y+z))\neq\emptyset \ \ \ \ {\rm for \ all} \ x,y\in R,\\
a+R=R+a=R \ \ \ \ {\rm for \ all} \ a\in R;
\end{array}
\]
 \item[(ii)] \ $(R,\cdot )$ is an $H_v$-semigroup;
 \item[(iii)] \ ``$\cdot$'' is {\em weak distributive} with respect to ``$+$'', i.e.,
for all $x,y,z\in R$:
\[
\begin{array}{l}
( x\cdot (y+z))\cap (x\cdot y+x\cdot z)\neq\emptyset ,\\
( (x+y)\cdot z)\cap (x\cdot z+y\cdot z)\neq\emptyset .
\end{array}
\]
\end{itemize}
\end{definition}
\begin{definition}(cf. Vougiouklis \cite{20}). A non-empty set $M$ is called
an {\em $H_v$-module} over an $H_v$-ring $R$ if $(M,+)$ is a weak
commutative $H_v$-group and there exists a map
 \[\cdot : R \times M \longrightarrow {\cal P}^*(M), \,
   (r,x) \mapsto  r \cdot x  \]
such that for all $a,b \in R$ and $x,y \in M$, we have
\[
\begin{array}{c}
(a \cdot (x+y))\cap (a \cdot x+a \cdot y)\not = \emptyset , \\
((a+b)\cdot x) \cap (a \cdot x+b \cdot x) \not = \emptyset , \\
((ab)\cdot x) \cap (a \cdot (b \cdot x)) \not = \emptyset .
\end{array}
\]
\end{definition}

We note that an $H_v$-module is a generalization of a module. For
more definitions, results and applications on $H_v$-modules, we
refer the reader to \cite{9, 18, 20}. Note that by using fuzzy
sets, we can consider the structure of $H_v$-module on any
ordinary module.

\begin{example}\rm (cf. Davvaz \cite{6}). Let $M$ be an
ordinary module over an ordinary ring $R$, and let $\mu_A$ be a
fuzzy set in $M$ and $\mu_B$ be a fuzzy set in $R$. We define
hyperoperations $\circ , \ *, \ \oplus $ and $\odot$ as follows:
\begin{itemize}
\item[ ] $a \circ b =\{t \in R \ | \ \mu_B(t)=\mu_B(a+b)\}$ for
         all $a,b \in R$,
\item[ ] $a*b=\{t \in R \ | \ \mu_B(t)=\mu_B(ab)
            \}$ for all $a,b \in R$,
\item[ ] $x\oplus y=\{s \in M \ | \
          \mu_A(s)=\mu_A(x+y)\}$ for all $x,y \in M$,
\item[ ] $r \odot
          x=\{s \in M \ | \ \mu_A(s)=\mu_A(r \cdot x)\}$ for all $r \in R$
and $x \in M$,
\end{itemize}
respectively.  Then
\begin{itemize}
 \item[(i)] $(R, \circ , *)$ is an $H_v$-ring.
 \item[(ii)] $(M,\oplus , \odot )$ is an $H_v$-module over the
 $H_v$-ring $(R, \circ ,*)$.
\end{itemize}
\end{example}

\begin{definition} Let $M$ be an $H_v$-module over an $H_v$-ring $R$. A non-empty
subset $S$ of $M$ is called an {\em $H_v$-submodule} of $M$ if the
following axioms hold:
\begin{itemize}
\item[(i)] \ $(S,+)\,$ is an $H_v$-subgroup of $(M,+)$,
\item[(ii)]\ $R\cdot S\subseteq S$.
\end{itemize}
\end{definition}

\begin{definition} Let $M_1$ and $M_2$ be two $H_v$-modules over
an $H_v$-ring $R$. A mapping $f$ from $M_1$ into $M_2$ is called a
{\em homomorphism} if for all $x,y \in M_1$ and $r \in R$,
\[ f(x+y)=f(x)+f(y) \ \ {\rm and} \ \ f(r\cdot x)=r \cdot f(x). \]
The homomorphism $f$ is said to be {\it strong on the left} if
 \[ f(z)\in f(x)+f(y) \ \Longrightarrow \ \exists x' \in M_1 \ : \
     f(x)=f(x') \ {\rm and} \ z \in x'+y. \]
Similarly, we can define a homomorphism which is {\it strong on
the right}. If a homomorphism $f$ is strong on the right and
left, we say $f$ is a {\it strong homomorphism}.
\end{definition}

\begin{proposition} Let $M_1$ and $M_2$ be two $H_v$-modules over
an $H_v$-ring $R$ and $f:M_1 \longrightarrow M_2$ a strong
epimorphism. If $N$ is an $H_v$-submodule of $M_2$, then
$f^{-1}(N)$ is an $H_v$-submodule of $M_1$.
\end{proposition}
\begin{proof} Assume that $x_1,x_2 \in f^{-1}(N)$. Then there
exists $y_1,y_2 \in N$ such that $f(x_1)=y_1$, $f(x_2)=y_2$ and so
$f(x_1+x_2)=y_1+y_2$. Hence for every $x \in x_1+x_2$ we have
$f(x) \in y_1+y_2 \subseteq N$ which implies $x \in f^{-1}(N)$ and
so $x_1+x_2 \subseteq f^{-1}(N)$. Therefore
    $x_1+f^{-1}(N) \subseteq f^{-1}(N)$ for every $x_1 \in
f^{-1}(N).$ Now, we show that $f^{-1}(N) \subseteq
x_1+f^{-1}(N)$. For $z \in f^{-1}(N)$, there exists $y \in N$ such
that $f(z)=y$. Since $y,y_1 \in N$ and $N$ is an $H_v$-subgroup of
$M_2$,  there exists $b \in N$ such that $y \in y_1+b$. Since $f$
is onto, there exists $a \in M_1$ such that $f(a)=b$ or $a \in
f^{-1}(b)$. Hence we have $f(z)\in f(x_1)+f(a)$. Since $f$ is a
strong homomorphism,
 there exists $a' \in M_1$ such that $f(a)=f(a')$ and $z \in
x_1+a'$. Since $f(a')=b \in N$, we have $a' \in f^{-1}(N)$ and $z
\in x_1+f^{-1}(N)$, and so $f^{-1}(N)\subseteq x_1+f^{-1}(N)$.
Therefore we have $f^{-1}(N)=x_1+f^{-1}(N)$. Similarly, we obtain
$f^{-1}(N)=f^{-1}(N)+x_1$. Thus the condition (i) of Definition
3.4 is satisfied.\\
\indent For the condition (ii), let $r \in R$ and $x \in
f^{-1}(N)$, then $f(x) \in N$, and so $r \cdot f(x) \subseteq N$
or $f(r\cdot x)\subseteq N$, which implies $r \cdot x \subseteq
f^{-1}(N)$. Therefore the condition (ii) of Definition 3.4 is
satisfied.
\end{proof}

In \cite{7}, Davvaz applied the concept of fuzzy sets to the
algebraic hyperstructures. In particular, he defined the concept
of a fuzzy $H_v$-submodule of an $H_v$-module which is a
generalization of the concept of a fuzzy submodule (see \cite{6}),
and he studied further properties  in \cite{8}, \cite{9} and
\cite{10}.

\begin{definition}(cf. Davvaz \cite{6}). Let $M$ be an $H_v$-module over
an  $H_v$-ring $R$ and $\mu$ a fuzzy set in $M$. Then $\mu$ is
said to be a {\em fuzzy $H_v$-submodule} of $M$ if the following
axioms hold:
\begin{itemize}
 \item[(1)] $\min\{\mu(x),\mu(y)\}\leq\inf\{\mu(z) \ | \ z\in
     x+y\}$ \ \ \ \ for all $x,y\in M$,
 \item[(2)] for all $x,a\in M$
there exists $y\in M$ such that $x\in a+y$ and
\[
\min\{\mu(a),\mu(x)\}\leq\mu(y),
\]
\item[(3)] for all $x,a\in M$ there exists $z\in M$ such that
$x\in z+a$ and
\[
\min\{\mu(a),\mu(x)\}\leq\mu (z),
\]
\item[(4)] $\mu(x)\leq\inf\{\mu(z) \ | \ z\in r\cdot x\}$ for all
$x\in M$ and $r\in R$.
\end{itemize}
\end{definition}

\section{Intuitionistic fuzzy $H_v$-submodules} 

In what follows, let $M$ denote an $H_v$-module over an $H_v$-ring
$R$ unless otherwise specified. We start by defining the notion of
intuitionistic fuzzy $H_v$-submodules.

\begin{definition}\label{def31} An intuitionistic fuzzy set $A=(\mu_A,
\lambda_A)$ in $M$ is called an {\em intuitionistic fuzzy
$H_v$-submodule} of $M$ if
\begin{itemize}
 \item[(1)] $\min\{\mu_A (x),\mu_A(y)\}\leq\inf\{\mu_A(z)\ |\ z\in
x+y\}$ \ for all $x,y\in M$,
 \item[(2)] for all $x,a\in M$ there
exist $y,z\in M$ such that\ $x\in (a+y)\cap (z+a)$\ and
\[
\min\{\mu_A(a),\mu_A(x)\}\leq\min\{\mu_A(y),\mu_A(z)\},
\]
 \item[(3)] $\mu_A(x)\leq\inf\{\mu_A(z)\ |\ z\in r\cdot x\}$ for
all $x \in M$ and $r\in R$,
 \item[(4)] $\sup\{\lambda_A(z)\ |\ z\in
x+y\}\leq\max\{\lambda_A(x),\lambda_A(y)\}$\ for all $x,y\in M$,
\item[(5)] for all $x,a\in M$ there exist $y,z\in M$ such that\
$x\in (a+y)\cap (z+a)$\ and
\[
\max\{\lambda_A(y),\lambda_A(z)\}\leq\max\{\lambda_A(a),
\lambda_A(x)\},
\]
\item[(6)] $\sup\{\lambda_A(z)\ |\ z\in r\cdot
x\}\leq\lambda_A(x)$ for all $x \in M$ and $r\in R$.
\end{itemize}
\end{definition}

\begin{lemma} If $A=(\mu_A,\lambda_A)$ is an
intuitionistic fuzzy $H_v$-submodule of $M$, then so is \ $\Box
A=(\mu_A,\mu^c_A)$.
\end{lemma}

\begin{proof} It is sufficient to show that $\mu_A^c$ satisfies the
conditions (4),(5), and (6) of Definition 4.1. For $x,y\in M$ we
have
$$
\min\{\mu_A(x),\mu_A(y)\}\leq\inf\{\mu_A(z)\ |\ z\in x+y\}
$$
and so $ \min\{1-\mu_A^c(x),1-\mu_A^c(y)\}\leq\inf\{1-\mu_A^c(z)\
|\ z\in x+y\}.
$
Hence
$$
\min\{1-\mu_A^c(x),1-\mu_A^c(y)\}\leq 1-\sup\{\mu_A^c(z)\ |\ z\in
x+y\}
$$
which implies
$
\sup\{\mu_A^c(z)\ |\ z\in x+y\}\leq 1-\min\{1-\mu_A^c(x),
1-\mu_A^c (y)\}.
$
Therefore
$$
\sup\{\mu_A^c(z)\ |\ z\in x+y\}\leq\max\{\mu_A^c(x),\mu_A^c(y)\},
$$
and thus the condition (4) of Definition 4.1 is valid.

Now, let $a,x\in M.$ Then there exist $y,z\in M$ such that $x\in
(a+y)\cap (z+a)$ and $
\min\{\mu_A(a),\mu_A(x)\}\leq\min\{\mu_A(y),\mu_A(z)\}.$ It
follows that
\[
\min\{1-\mu_A^c(a),1-\mu_A^c(x)\}\leq\min\{1-\mu_A^c(y),1-\mu_A^c(z)\}
\]
so that
\[
\max\{\mu_A^c(y),\mu_A^c(z)\}\leq\max\{\mu_A^c(a),\mu_A^c(x)\}.
\]
Hence the condition (5) of Definition 4.1 is satisfied.

For the condition (6), let $x\in M$ and $r\in R$. Since $\mu_A$ is
a fuzzy $H_v$-submodule of $M$, we have
\[
\mu_A(x)\leq\inf\{\mu_A(z)\ |\ z\in r\cdot x\}
\]
and so
\[
1-\mu_A^c(x)\leq\inf\{1-\mu_A^c(z)\ |\ z\in r\cdot x\},
\]
which implies
\[
\sup\{\mu_A^c(z)\ |\ z\in r\cdot x\}\leq\mu_A^c(x).
\]
Therefore the condition (6) of Definition 4.1 is satisfied.
\end{proof}

\begin{lemma} If $A=(\mu_A,\lambda_A)$ is an
intuitionistic fuzzy $H_v$-submodule of $M$, then so is \
$\Diamond A=(\lambda^c_A,\lambda_A)$.
\end{lemma}

\begin{proof} The proof is similar to the proof of Lemma 4.2.
 \end{proof}

Combining the above two lemmas it is not difficult to verify that
the following theorem is valid.

\begin{theorem} $A=(\mu_A, \lambda_A)$ is an
intuitionistic fuzzy $H_v$-submodule of $M$ if and only if \ $\Box
A$ and $\Diamond A$ are intuitionistic fuzzy $H_v$-submodules of
$M$. \hfill $\Box$
\end{theorem}
\begin{corollary} $A=(\mu_A,\lambda_A)$ is an
intuitionistic fuzzy $H_v$-submodule of $M$ if and only if \
$\mu_A$ and $\lambda^c_A$ are fuzzy $H_v$-submodules of $M$.
\end{corollary}

\begin{definition} For any $t\in [0,1]$ and fuzzy set
$\mu$ in $M$, the set
\[
U(\mu ;t)=\{x\in M\ |\ \mu(x)\geq t\} \ \ ({\rm resp. } \ L(\mu
;t) =\{x\in M\ |\ \mu(x)\leq t\} )
\]
is called an {\em upper} (resp. {\em lower}) {\em $t$-level cut}
of $\mu$.
\end{definition}

\begin{theorem} If $A=(\mu_A,\lambda_A)$ is an
intuitionistic fuzzy $H_v$-submodule of $M$, then  the sets
$U(\mu_A ;t)$ and $L(\lambda_A;t)$ are $H_v$-submodules of $M$ for
every $t\in Im(\mu_A)\cap Im(\lambda_A).$
\end{theorem}

\begin{proof} Let $t\in Im(\mu_A)\cap Im(\lambda_A)\subseteq [0,1]$
and let $x,y\in U(\mu_A ;t)$. Then $\mu_A(x)\geq t$ and
$\mu_A(y)\geq t$ and so $\min\{\mu_A (x),\mu_A(y) \}\geq t$. It
follows from the condition (1) of Definition 4.1 that
$\inf\{\mu_A(z)\ |\ z\in x+y\}\geq t$. Therefore $z\in U(\mu_A;t)$
for all $z\in x+y,$ and so $x+y\subseteq U(\mu_A;t)$. Hence
$a+U(\mu_A;t)\subseteq U(\mu_A;t)$ and $U(\mu_A;t)+a\subseteq
U(\mu_A;t)$ for all $a\in U(\mu_A;t).$ Now, let $x\in U(\mu_A;t).$
Then there exist $y,z\in M$ such that $x\in (a+y)\cap (z+a)$ and
$\min\{\mu_A(x),\mu_A(a)\}\leq\min\{\mu_A(y),\mu_A(z)\}$. Since
$x,a\in U(\mu_A;t)$, we have $t\leq\min\{\mu_A(x),\mu_A(a)\}$ and
so $t\leq\min\{\mu_A(y),\mu_A(z)\} ,$ which implies $y\in
U(\mu_A;t)$ and $z\in U(\mu_A;t)$. This proves that
$U(\mu_A;t)\subseteq a+U(\mu_A;t)$ and $U(\mu_A;t)\subseteq
U(\mu_A;t)+a$.

Now, for every $r\in R$ and $x\in U(\mu_A;t)$ we show that $r\cdot
x\subseteq U(\mu_A;t)$. Since $A$ is an intuitionistic fuzzy
$H_v$-submodule of $M$, we have
\[
t\leq \mu_A(x)\leq\inf\{\mu_A(z)\ |\ z\in r\cdot x\}.
\]
Therefore, for every $z\in r\cdot x$ we get $\mu_A(z)\geq t$ which
implies $z\in U(\mu_A;t)$, so $r\cdot x\subseteq U(\mu_A;t)$.

If $x,y\in L(\lambda_A;t)$, then
$\max\{\lambda_A(x),\lambda_A(y)\}\leq t$. It follows from the
condition (4) of Definition 4.1 that $\sup\{\lambda_A(z)\ |\ z\in
x+y\}\leq t$. Therefore for all $z\in x+y$ we have $z\in
L(\lambda_A;t)$, so $x+y\subseteq L(\lambda_A;t)$. Hence for all
$a\in L(\lambda_A;t)$ we have $a+L(\lambda_A;t)\subseteq
L(\lambda_A;t)$ and $L(\lambda_A;t)+a \subseteq L(\lambda_A;t)$.
Now, let $x\in L(\lambda_A;t).$ Then there exist $y,z\in M$ such
that $x\in (a+y)\cap (z+a)$ and
$\max\{\lambda_A(y),\lambda_A(z)\}\leq\max\{\lambda_A(a),\lambda_A(x)\}$.
Since $x,a\in L(\lambda_A;t)$, we have
$\max\{\lambda_A(a),\lambda_A(x)\}\leq t$ and so
$\max\{\lambda_A(y),\lambda_A(z)\}\leq t$. Thus $y\in
L(\lambda_A;t)$ and $z\in L(\lambda_A;t)$. Hence
$L(\lambda_A;t)\subseteq a+L(\lambda_A;t)$ and
$L(\lambda_A;t)\subseteq L(\lambda_A;t)+a$.

Now, we show that $r\cdot x\subseteq L(\lambda_A;t)$ for every
$r\in R$ and $x\in L(\lambda_A;t)$. Since $A$ is an intuitionistic
fuzzy $H_v$-submodule of $M$, we have
\[
\sup\{\lambda_A(z)\ |\ z\in r\cdot x\}\leq\lambda_A(x)\leq t.
\]
Therefore, for every $z\in r\cdot x$ we get $\lambda_A(z)\leq t,$
which implies $z\in L(\lambda_A;t)$, so $r\cdot x\subseteq
L(\lambda_A;t)$.
\end{proof}

\begin{theorem} If $A=(\mu_A,\lambda_A)$ is an
intuitionistic fuzzy set in $M$ such that all non-empty level sets
$U(\mu_A;t)$ and $L(\lambda_A;t)$ are $H_v$-submodules of $M$,
then $A=(\mu_A,\lambda_A)$ is an intuitionistic fuzzy
$H_v$-submodule of $M$.
\end{theorem}

\begin{proof} Assume that all non-empty level sets $U(\mu_A;t)$
and $L(\lambda_A;t)$ are $H_v$-submodules of $M$. If
$t_0=\min\{\mu_A(x),\mu_A(y)\}$ and
$t_1=\max\{\lambda_A(x),\lambda_A(y)\}$ for $x,y\in M$, then
$x,y\in U(\mu_A;t_0)$ and $x,y\in L(\lambda_A;t_1)$. So
$x+y\subseteq U(\mu_A;t_0)$ and $x+y\subseteq L(\lambda_A;t_1)$.
Therefore for all $z\in x+y$ we have $\mu_A(z)\geq t_0$ and
$\lambda_A(z)\leq t_1$, i.e.,
\[
\inf\{\mu_A(z)\ |\ z\in x+y\}\geq\min\{\mu_A(x),\mu_A(y)\}
\]
and
\[
\sup\{\lambda_A(z)\ |\ z\in
x+y\}\leq\max\{\lambda_A(x),\lambda_A(y)\},
\]
which verify the conditions (1) and (4) of Definition 4.1.

Now, if $t_2=\min\{\mu_A(a),\mu_A(x)\}$  for $x,a\in M$, then
$a,x\in U(\mu_A;t_2)$. So there exist $y_1,z_1\in U(\mu_A;t_2)$
such that $x\in a+y_1$ and $x\in z_1+a$. Also we have
$t_2\leq\min\{\mu_A(y_1),\mu_A(z_1)\}$. Therefore the condition
(2) of Definition 4.1 is verified. If we put
$t_3=\max\{\lambda_A(a),\lambda_A(x)\}$ then $a,x\in
L(\lambda_A;t_3)$. So there exist $y_2,z_2\in L(\lambda_A;t_3)$
such that $x\in a+y_2$ and $x\in z_2+a$ and we have
$\max\{\lambda_A(y_2),\lambda_A(y_2)\}\leq t_3$, and so the
condition (5) of Definition 4.1 is verified.

Now, we verify the conditions (3) and (6). Let $t_4=\mu_A(x)$ and
$t_5=\lambda_A(x)$ for some $x\in M$ and let $r \in R$. Then $x\in
U(\mu_A;t_4)$ and $x\in L(\lambda_A, t_5)$. Since $U(\mu_A;t_4)$
and $L(\lambda_A,t_5)$ are $H_v$-submodules of $M$, we get $r\cdot
x \subseteq U(\mu_A;t_4)$ and $r\cdot x\subseteq
L(\lambda_A,t_5)$. Therefore for every $z\in r\cdot x$ we have
$z\in U(\mu_A;t_4)$ and $z\in L(\lambda_A,t_5)$ which imply
$\mu_A(z)\geq t_4$ and $\lambda_A(z)\leq t_5$. Hence
\[
\inf\{\mu_A(z)\ |\ z\in r\cdot x\}\geq t_4=\mu_A(x)
\]
and
\[
\sup\{\lambda_A(z)\ |\ z\in r\cdot x\}\leq t_5=\lambda_A(x).
\]
This completes the proof.
\end{proof}

\begin{corollary} Let $S$ be an $H_v$-submodule of
an $H_v$-module $M$. If fuzzy sets $\mu$ and $\lambda$ in $M$ are
defined by
\[
\mu(x)=\left\{\begin{array}{cl}\alpha_0 &\text{if \, $x\in S$},\\
               \alpha_1 &\text{if \, $x\in M\setminus S$,}\end{array}
\right. \qquad\lambda(x)=\left\{\begin{array}{cl}\beta_0
&\text{if \, $x\in S$},\\
                  \beta_1 &\text{if \, $x\in M\setminus S$,}\end{array}
\right.
\]
where $\,0\leq\alpha_1< \alpha_0$, \ $0\leq\beta_0<\beta_1\,$ and
$\,\alpha_i+\beta_i\leq 1$ for $i=0,1$. Then $A=(\mu,\lambda )$ is
an intuitionistic fuzzy $H_v$-submodule of $M$ and
$\,U(\mu;\alpha_0 )=S=L(\lambda ;\beta_0)$.
\end{corollary}

\begin{corollary}
Let $\chi_{_S}$ be the characteristic function of an
$H_v$-submodule $S$ of $M$. Then $A=(\chi_{_S},\chi^c_{_S})$ is an
intuitionistic fuzzy $H_v$-submodule of $M$.  
\end{corollary}

\begin{theorem} If $\,A=(\mu_{A},\lambda_A)\,$ is an
intuitionistic fuzzy $H_v$-submodule of $M$, then
\\[2mm]
\centerline{$\mu_{A}(x)=\sup\{\alpha \in [0,1]\ |\ x\in
U(\mu_{A};\alpha )\}$}
\\
and\\
\centerline{$ \lambda_A(x)=\inf\{\alpha \in [0,1]\ |\ x\in
L(\lambda_A ;\alpha)\}$ } for all $\,x\in M.$
\end{theorem}

\begin{proof} Let $\,\delta=\sup\{\alpha\in [0,1]\ |\ x\in
U(\mu_{A};\alpha )\}\,$ and let $\,\varepsilon >0\,$ be given.
Then $\delta -\varepsilon <\alpha$ for some $\,\alpha \in [0,1]$
such that $\,x\in U(\mu_{A};\alpha)$. This means that $\,\delta
-\varepsilon <\mu_{A}(x)\,$ so that $\,\delta \leq\mu_{A}(x)\,$
since $\,\varepsilon\,$ is arbitrary.

We now show that $\,\mu_{A}(x)\leq\delta.$ If
$\,\mu_{A}(x)=\beta$, then $\,x\in U(\mu_{A};\beta)\,$ and so
\[
\beta\in\{\alpha\in [0,1]\ |\ x\in U(\mu_{A};\alpha )\}.
\]
Hence
\[
\mu_{A}(x)=\beta\leq\sup\{\alpha\in [0,1]\ |\ x\in
U(\mu_{A};\alpha)\}=\delta .
\]
Therefore
\[
\mu_{A}(x)=\delta=\sup\{\alpha\in [0,1]\ |\ x\in U(\mu_{A}
;\alpha)\}.
\]

Now let $\,\eta =\inf\{\alpha \in [0,1]\ |\ x\in
L(\lambda_A;\alpha)\}$. Then
\[
\inf\{\alpha\in [0,1]\ |\ x\in L(\lambda_A;\alpha)\}<\eta
+\varepsilon
\]
for any  $\,\varepsilon >0,\,$ and so $\,\alpha <\eta
+\varepsilon\,$ for some $\,\alpha\in [0,1]\,$ with $\,x\in
L(\lambda_A ;\alpha)$. Since $\lambda_A(x)\leq\alpha\,$ and
$\,\varepsilon\,$ is arbitrary, it follows that
$\,\lambda_A(x)\leq\eta$.

To prove $\,\lambda_A(x)\geq\eta$, let $\,\lambda_A(x)=\zeta$.
Then $\,x\in L(\lambda_A ;\zeta)\,$ and thus
$\,\zeta\in\{\alpha\in [0,1]\ |\ x\in L(\lambda_A ;\alpha)\}$.
Hence
\[
\inf\{\alpha\in [0,1]\ |\ x\in L(\lambda_A ;\alpha)\}\leq\zeta,
\]
i.e. $\;\eta\leq\zeta =\lambda_A(x).$ Consequently
\[
\lambda_A(x)=\eta=\inf\{\alpha\in [0,1]\ |\ x\in L(\lambda_A
;\alpha)\},
\]
which completes the proof.
\end{proof}

\begin{definition} A fuzzy set $\mu$ in a set $X$ is said to have
 {\it sup property} if for every non-empty subset $S$ of $X$,
there exists $x_0 \in S$ such that
\[
\mu (x_0)=\displaystyle \sup_{x \in S} \{ \mu (x)\}.
\]
\end{definition}

\begin{proposition} Let $M_1$ and $M_2$ be two $H_v$-modules over
an $H_v$-ring $R$ and $f: M_1 \longrightarrow M_2$ be a
surjection. If $A=(\mu_A, \lambda_A )$ is an intuitionistic fuzzy
$H_v$-submodule of $M_1$ such that $\mu_A$ and $\lambda_A$ have
sup property, then
\begin{itemize}
 \item[\rm (i)] $f(U(\mu_A;t))=U(f(\mu_A);t),$
 \item[\rm (ii)] $f(L(\lambda_A;t))\subseteq L(f(\lambda_A);t)$
\end{itemize}
\end{proposition}

\begin{proof} (i) We have
\[
\begin{array}{ll}
y \in U(f(\mu_A);t) & \Longleftrightarrow f(\mu_A)(y)\geq t \\
     & \Longleftrightarrow \displaystyle \sup_{x \in f^{-1}(y)}\{ \mu_A(x)\} \geq t \\
     & \Longleftrightarrow \exists x_0\in f^{-1}(y), \ \mu_A (x_0)\geq t \\
     & \Longleftrightarrow \exists x_0 \in f^{-1}(y), \ x_0 \in U(\mu_A;t)\\
     & \Longleftrightarrow f(x_0)=y, \ x_0 \in U(\mu_A;t) \\
     & \Longleftrightarrow y \in f(U(\mu_A;t)).
\end{array} \]
(ii) We have
\[
\begin{array}{ll}
y \in L(f(\lambda_A);t) & \Longrightarrow f(\lambda_A)(y)\leq t \\
     & \Longrightarrow \displaystyle \sup_{x \in f^{-1}(y)}\{ \lambda_A(x)\} \leq t \\
     & \Longrightarrow  \lambda_A (x)\leq t\ {\rm for \ all} \ x \in f^{-1}(y) \\
     & \Longrightarrow x \in L(\lambda_A;t)\ {\rm for \ all} \ x \in f^{-1}(y)\\
     & \Longrightarrow y \in f(L(\lambda_A;t)).
\end{array}
\]
\end{proof}

\begin{proposition} Let $M_1$ and $M_2$ be two $H_v$-modules over
an $H_v$-ring $R$ and $f: M_1 \longrightarrow M_2$ be a map. If
$B=(\mu_B, \lambda_B )$ is an intuitionistic fuzzy $H_v$-submodule
of $M_2$, then
\begin{itemize}
 \item[\rm (i)] $f^{-1}(U(\mu_B;t))=U(f^{-1}(\mu_B);t),$
 \item[\rm (ii)] $f^{-1}(L(\lambda_B;t))=L(f^{-1}(\lambda_B);t)$
\end{itemize}
for every $t\in [0,1].$
\end{proposition}

 \begin{proof} (i) We have
\[
\begin{array}{ll}
x \in U(f^{-1}(\mu_B);t) & \Longleftrightarrow f^{-1}(\mu_B)(x)\geq t \\
     & \Longleftrightarrow \mu_B(f(x)) \geq t \\
     & \Longleftrightarrow  f(x)\in U(\mu_B;t)\\
     & \Longleftrightarrow x \in f^{-1}(U(\mu_B;t)).
\end{array}
\]
(ii) We have
\[
\begin{array}{ll}
x \in L(f^{-1}(\lambda_B);t) & \Longleftrightarrow f^{-1}(\lambda_B)(x)\leq t \\
     & \Longleftrightarrow \lambda_B(f(x)) \leq t \\
     & \Longleftrightarrow  f(x)\in L(\lambda_B;t)\\
     & \Longleftrightarrow x \in f^{-1}(L(\lambda_B;t)).
\end{array}
\]
\end{proof}

\begin{definition} Let $f$ be a map from a set $X$ to a set $Y$.
If $B=(\mu_B, \lambda_B)$ is an intuitionistic fuzzy set in $Y,$
then the {\it inverse image} of $B$ under $f$ is defined by:
\[
f^{-1}(B)=(f^{-1}(\mu_B), f^{-1}(\lambda_B)).
\]
\end{definition}

It is easy to see that $f^{-1}(B)$ is an intuitionistic fuzzy set
in $X$.

\begin{corollary} Let $M_1$ and $M_2$ be two $H_v$-modules over an
$H_v$-ring $R$ and $f: M_1 \longrightarrow M_2$ be a strong
epimorphism. If $B=(\mu_B, \lambda_B )$ is an intuitionistic fuzzy
$H_v$-submodule of $M_2$, then $f^{-1}(B)$ is an intuitionistic
fuzzy $H_v$-submodule of $M_1$.
\end{corollary}

\begin{proof}
Assume that $B=(\mu_B, \lambda_B )$ is an intuitionistic fuzzy
$H_v$-submodule of $M_2$. By Theorem 4.7, we know that  the sets
$U(\mu_B;t)$ and $L(\lambda_B;t)$ are $H_v$-submodules of $M_2$
for every $t \in Im(\mu_B) \cap Im(\lambda_B).$ It follows from
Proposition 3.6 that $f^{-1}(U(\mu_B;t))$ and
$f^{-1}(L(\lambda_B;t))$ are $H_v$-submodules of $M_1$. Using
Proposition 4.14, we have
$$f^{-1}(U(\mu_B;t))=U(f^{-1}(\mu_B);t),$$
$$f^{-1}(L(\lambda_B;t))=L(f^{-1}(\lambda_B);t).$$
Now by Theorem 4.8, the proof is completed.
\end{proof}

\section{On fundamental modules} 

The main tools in the theory of $H_v$-structures are the
fundamental relations. Consider an $H_v$-module $M$ over an
$H_v$-ring $R$. If the relation $\gamma^*$ is the smallest
equivalence relation on $R$ such that the quotient $R/\gamma^*$,
the set of all equivalence classes, is a ring, we say that
$\gamma^*$ is the {\it fundamental equivalence relation} on $R$
and $R/\gamma^*$ is the {\it fundamental ring} (see \cite{17,
19}). The fundamental relation $\epsilon^*$ on $M$ over $R$ is the
smallest equivalence relation on $M$ such that $M/\epsilon^*$ is a
module over the ring $R/\gamma^*$. Let ${\cal U}$ be the set of
all expressions consisting of finite hyperoperations either on $R$
and $M$ or the external hyperoperation applied to finite sets of
elements of $R$ and $M$. We define the relation $\epsilon$ on $M$
as follows:
\[
a \, \epsilon \, b \ {\rm if \ and \ only \ if } \ \{ a,b\}
\subseteq u \ \ {\rm for \ some } \ u\in {\cal U}.
\]
Let us denote  $\widehat {\epsilon}$  the transitive closure of
$\epsilon$. Then we can rewrite the definition of
$\widehat{\epsilon}$ on $M$ as follows:\\
\indent  $a \, \widehat{\epsilon} \, b$ if and only if there exist
$z_1, \ldots , z_{n+1}\in M$ with $z_1=a$, $z_{n+1}=b$ and $ u_1,
\ldots , u_n\in {\cal U}$ such that
\[
\{z_i, z_{i+1}\} \subseteq u_i \ \ \ \ (i=1, \ldots ,n).
\]

\begin{theorem} {\rm ( cf. Vougiouklis \cite{20})}. The fundamental
relation $\epsilon^*$ is the transitive closure of the relation
$\epsilon$.
\end{theorem}

Suppose $\gamma^* (r)$ is the equivalence class containing $r\in
R,$ and $\epsilon^*(x)$ the equivalence class containing $x \in
M$. On $M/\epsilon^*$, the sum $\oplus $ and the external product
$\odot $ using the $\gamma^*$ classes in $R$ are defined as
follows:
\[
\epsilon^* (x) \oplus \epsilon^* (y)=\epsilon^*(c) \ {\rm for \
all } \  c\in \epsilon^* (x)+\epsilon^*(y),
\]
\[
\gamma^* (r) \odot \epsilon^* (x)=\epsilon^*(d) \  {\rm for \ all
} \  d\in \gamma^* (r)\cdot \epsilon^*(x).
\]
The kernel of the canonical map $\varphi : M \longrightarrow
M/\epsilon^*$ is called the {\it core} of $M$ and is denoted by
$\omega_M$. Here we also denote $\omega_M$ the zero element of
$M/\epsilon$. We have
\[
\omega_M=\epsilon^*(0), \ \ {\rm and} \
\epsilon^*(-x)=-\epsilon^*(x) \ \ {\rm for \ all} \ x \in M.
\]

\begin{definition} Let $M$ be an $H_v$-module over an $H_v$-ring $R$
and let $A=(\mu_A, \lambda_A)$ be an intuitionistic fuzzy
$H_v$-submodule of $M$. The intuitionistic fuzzy set
$A/\epsilon^*=(\overline {\mu_A}^{\epsilon^*}, \ \underline
{\lambda_A}_{\epsilon^*})$ is defined as follows:
\[
\overline {\mu_A}^{\epsilon^*} :M/\epsilon^* \longrightarrow [0,1]
\]
\[
\overline {\mu_A}^{\epsilon^*} (\epsilon^*(x))=\left \{
\begin{array}{ll} \displaystyle \sup_{a \in \epsilon^*(x)} \{\mu_A
(a)\} & {\rm if}
\ \epsilon^*(x)\not = \omega_M \\
1 & {\rm if} \ \epsilon^*(x)=\omega_M
\end{array}
\right.
\]
and
\[
\underline {\lambda_A}_{\epsilon^*} :M/\epsilon^* \longrightarrow
[0,1]
\]
\[
\underline {\lambda_A}_{\epsilon^*}(\epsilon^*(x))=\left \{
\begin{array}{ll} \displaystyle \inf_{a \in \epsilon^*(x)}
\{\lambda_A (a)\} & {\rm if}
\ \epsilon^*(x)\not = \omega_M \\
0 & {\rm if} \ \epsilon^*(x)=\omega_M.
\end{array}
\right.
\]
\end{definition}

In the following we show that
\[
0\leq \overline {\mu_A}^{\epsilon^*}(\epsilon^*(x))+ \underline
{\lambda_A}_{\epsilon^*}(\epsilon^*(x)) \leq 1,
\]
for all $\epsilon^*(x)\in M/\epsilon^*$. \\
\indent If $\epsilon^*(x)=\omega_M$, then the above inequalities
are clear. Assume that $x \in H$ and $\epsilon^*(x)\not
=\omega_M$. Since $0\leq \mu_A(a)$ and $0\leq \lambda_A(a)$ for
all $a \in \epsilon^*(x)$, we have
\[
0\leq \displaystyle \sup_{a \in \epsilon^*(x)}\{ \mu_A(a) \} +
\displaystyle \inf_{a \in \epsilon^*(x)}\{ \lambda_A(a) \}
\]
or
\[
0\leq \overline {\mu_A}^{\epsilon^*}(\epsilon^*(x))+ \underline
{\lambda_A}_{\epsilon^*}(\epsilon^*(x)).
\]

On the other hand,  we have
\[
\mu_A(a)+\lambda_A(a) \leq 1 \ \ \ {\rm or} \ \ \ \mu_A(a) \leq
1-\lambda_A(a), \] for all $a \in \epsilon^*(x),$ and so
\[
\begin{array}{ll}
\overline {\mu_A}^{\epsilon^*}(\epsilon^*(x)) & =\displaystyle
\sup_{a \in \epsilon^*(x)}\{ \mu_A(a) \} \\
 & \leq \displaystyle \sup_{a \in \epsilon^*(x)}\{1-\lambda_A(a)
 \}\\
& =1- \displaystyle \inf_{a \in \epsilon^*(x)}\{\lambda_A(a)
 \}\\
& =1- \underline {\lambda_A}_{\epsilon^*}(\epsilon^*(x)).
\end{array}
\]
Hence $\overline {\mu_A}^{\epsilon^*}(\epsilon^*(x))+ \underline
{\lambda_A}_{\epsilon^*}(\epsilon^*(x)) \leq 1$.

\begin{theorem}{\rm (cf. Davvaz \cite{6}).} Let $M$ be an
$H_v$-module over an $H_v$-ring $R$ and let $\mu$ be a fuzzy
$H_v$-submodule of $M$. Then $\overline {\mu_A}^{\epsilon^*}$ is a
fuzzy submodule of the module $M/\epsilon^*$.
\end{theorem}

\begin{lemma} We have
\[
(\overline {\lambda_A^c}^{\epsilon^*})^c= \underline
{\lambda_A}_{\epsilon^*}.
\]
\end{lemma}
\begin{proof} If $\epsilon^*(x)=\omega_M$, then
\[
(\overline {\lambda_A^c}^{\epsilon^*})^c (\omega_H)= 1- (\overline
{\lambda_A^c}^{\epsilon^*}) (\omega_M)=0= \underline
{\lambda_A}_{\epsilon^*}(\omega_M).
\]
Now, assume that $\epsilon^*(x)\not =\omega_M.$ Then
\[
\begin{array}{ll}
( \overline {\lambda_A^c}^{\epsilon^*})^c (\epsilon^*(x)) & = 1-
(\overline {\lambda_A^c}^{\epsilon^*}) (\epsilon^*(x))\\
& =1- \displaystyle \sup_{a \in \epsilon^*(x)}\{\lambda_A^c(a)\}
\\
& =1- \displaystyle \sup_{a \in \epsilon^*(x)}\{1- \lambda_A(a)\}
\\
& = \displaystyle \inf_{a \in \epsilon^*(x)}\{\lambda_A(a)\} \\
& = \underline {\lambda_A}_{\epsilon^*}(\epsilon^*(x)).
\end{array}
\]
\end{proof}

\begin{theorem}
Let $M$ be an $H_v$-module over an $H_v$-ring $R$ and let
$A=(\mu_A, \lambda_A)$ be an intuitionistic fuzzy $H_v$-submodule
of $M$. Then $A/\epsilon^*=$ $(\overline {\mu_A}^{\epsilon^*},$
$\underline {\lambda_A}_{\epsilon^*})$ is an intuitionistic fuzzy
submodule of the fundamental module $M/\epsilon^*.$
\end{theorem}

\begin{proof} Suppose that $A=(\mu_A, \lambda_A)$ is an
intuitionistic fuzzy $H_v$-submodule of $M$. Using Lemma 4.3,
$\lambda_A^c$ is a fuzzy $H_v$-submodule of $M$ and by Theorem
5.3, $\overline {\mu_A}^{\epsilon^*}$ and $\overline
{\lambda_A^c}^{\epsilon^*}$ are fuzzy $H_v$-submodules of
$M/\epsilon^*$, and so $(\overline {\lambda_A^c}^{\epsilon^*})^c$
satisfies the conditions (4), (5), (6) of Definition 2.5. Hence by
Lemma 5.4, $\underline {\lambda_A}_{\epsilon^*}$ satisfies the
conditions (4),(5),(6) of Definition 2.5. Therefore
$A/\epsilon^*=(\overline {\mu_A}^{\epsilon^*}, \ \underline
{\lambda_A}_{\epsilon^*})$ is an intuitionistic fuzzy submodule of
$M/\epsilon^*$.
\end{proof}

\section{Conclusions}

As a generalization of fuzzy sets, the notion of intuitionistic
fuzzy sets was introduced by Atanassov \cite{1}, and applications
of intuitionistic fuzzy concepts have already been done by
Atanassov and others in algebra, topological space, knowledge
engineering, natural language, and neural network etc. Biswas
\cite{3} have applied the concept of intuitionistic fuzzy sets to
the theory of groups and studied intuitionistic fuzzy subgroups of
a group. The notion of an intuitionistic fuzzy subquasigroup of a
guasigroup was discussed by Kim, Dudek and Jun \cite{12}. Also the
concept of intuitionistic fuzzy ideals of semigroups was
considered by Kim and Jun \cite{13}. The concept of hyperstructure
first was introduced by Marty \cite{14}. Vougiouklis \cite{19}, in
the fourth AHA congress (1990), introduced the notion of
$H_v$-structures. Recently, present authors \cite{11} have
discussed the intuitionistic fuzzification of the concept of
subhyperquasigroups in a hyperquasigroup. The aim of this paper is
to introduce the notion of an intuitionistic fuzzy $H_v$-submodule
of an $H_v$-module, and to investigate related properties.
Characterizations of intuitionistic fuzzy $H_v$-submodules are
given. Our future work will focus on studying the intuitionistic
fuzzy structure of $H_v$-nearring modules.

\section{Acknowledgements}

The authors are highly grateful to referees and Professor Witold
Pedrycz, Editor-in-Chief, for their valuable comments and
suggestions for improving the paper.

\end{document}